\documentclass[12pt]{article}
\usepackage{amssymb}
\usepackage{amsmath}

\begin{document}

\begin{titlepage}
\title{\bf Lifts of Structures on Product Manifolds}
\author{ Mehmet Tekkoyun \footnote{Corresponding author. E-mail address: tekkoyun@pau.edu.tr; Tel: +902582953616; Fax: +902582963535}\\
{\small Department of Mathematics, Pamukkale University,}\\
{\small 20070 Denizli, Turkey}}
\date{\today}
\maketitle
\begin{abstract}
This article presents the further steps of the previously done
studies taking into consideration the k-th order extensions of a
complex manifold. In the previous studies higher order vertical
and complete lifts of structures on the complex manifold were
introduced. Presently, k-th extended spaces of a product manifold
have been set
 and the higher order vertical, complete, complete- vertical
and horizontal lifts of geometric structures on the product
manifold have been presented.
\end{abstract}
\end{titlepage}

\section{\textbf{Introduction}}

Lifting theory that permits to extend the differentiable structures has also
an important role in differential geometry (see \cite{bowman, civelek,
tekkoyun1, tekkoyun2, tekkoyun3} and there in). In \cite{civelek}, the
structure of extended vector bundles has been obtained, especially the
canonical $k$-th order extended vector bundle $\pi ^{k}$ of a vector bundle $%
\pi =(E,\pi ,M)$. Extensions of a complex manifold were defined and also
higher order vertical, complete and horizontal lifts of complex functions,
vector fields, 1-forms and tensor fields of type (1,1) and type (0,2) on any
complex manifold to its extension spaces were studied in \cite{tekkoyun1,
tekkoyun2, tekkoyun3}.

The paper is structured as follows. In section 2, we recall extended complex
manifolds and higher order vertical and complete lifts of differential
elements on any complex manifold to its extension spaces and also extended
Kaehlerian manifolds. In section 3, using structures obtained in \cite%
{civelek, tekkoyun1, tekkoyun2, tekkoyun3} we define the $k$-th order
extension $^{k}N$ of any product manifold $N$ of dimension $2m+1.$ Also we
obtain the higher order vertical, complete, complete-vertical and horizontal
lifts of functions, vector fields and 1-forms on $N$ to $^{k}N.$ Then we
find the higher order vertical and complete lifts of tensor field of type
(1,1) on $N$ to $^{k}N$.

Throughout the paper, all objects are assumed to be differentiable of class $%
C^{\infty }$ and the sum is taken over repeated indices. Unless otherwise
stated it will be accepted $0\leq r\leq k,$ $1\leq i\leq m.$ Also, $v$ and $%
c $ will denote the vertical and complete lifts of geometric structures
either from $^{k-1}M$ to $^{k}M$ or from $^{k-1}N$ to $^{k}N.$ The symbol $%
\mathrm{C}_{j}^{r}$ called combination is the binomial coefficient $\binom{r%
}{j}.$

\section{\textbf{Preliminaries}}

In this section, we recall $k$-th order extension of a complex manifold and
higher order vertical and complete lifts of geometric elements on complex
manifold to its extension spaces, and also extended Kaehlerian manifolds
given in \cite{tekkoyun1, tekkoyun2}.

\subsection{\textbf{Extended Complex Manifolds}}

Let $M$ be 2m-real dimensional manifold and $^{k}M$ its $k$-th order
extended manifold$.$ A tensor field $J_{k}$ on $^{k}M$ is called an \textit{%
extended almost complex structure} on $^{k}M$ if $J_{k}$ is endomorphism of
the tangent space $T_{p}(^{k}M)$ such that $(J_{k})^{2}=-I$ at every point $%
p $ of $^{k}M.$ An extended manifold $^{k}M$ with extended almost complex
structure $J_{k}$ is called \textit{extended almost complex manifold}. If $%
k=0$, $J_{0}$ is called \textit{almost complex structure}\textbf{\ }and a
manifold $^{0}M=M$ with almost complex structure $J_{0}$ is said to be
\textit{almost complex manifold}.

Let $(x^{ri},\,y^{ri})$ be a real coordinate system on a neighborhood $^{k}U$
of any point $p$ of $^{k}M$. In this situation, we respectively define by $\{%
\frac{\partial }{\partial x^{ri}},\frac{\partial }{\partial y^{ri}}\}$ and $%
\{dx^{ri},dy^{ri}\}$ the natural bases over $\mathbf{R}$ of tangent space $%
T_{p}(^{k}M)$ and cotangent space $T_{p}^{\ast }(^{k}M)$ of $^{k}M.$ The
manifold $^{k}M$ is called \textit{extended complex manifold} if there
exists an open covering $\left\{ ^{k}U\right\} $ of $^{k}M$ satisfying the
following condition: \newline
There is a local coordinate system $(x^{ri},\,y^{ri})$ on each $^{k}U$ such
that for each point of $^{k}U,$%
\begin{equation}
J_{k}(\frac{\partial }{\partial x^{ri}})=\frac{\partial }{\partial y^{ri}}%
,J_{k}(\frac{\partial }{\partial y^{ri}})=-\frac{\partial }{\partial x^{ri}}.
\label{1.1}
\end{equation}%
If $k=0$, then the manifold $^{0}M=M$ with almost complex structure $J_{0}$
is said to be \textit{complex manifold}. Let $z^{ri}=x^{ri}+\mathbf{i}%
\,y^{ri},\mathbf{i}=\sqrt{-1}$ be an extended complex local coordinate
system on a neighborhood $^{k}U$ of any point $p$ of $^{k}M.$ Then one puts:
$\frac{\partial }{\partial z^{ri}}=\frac{1}{2}\left\{ \frac{\partial }{%
\partial x^{ri}}-{\mathbf{i}}\frac{\partial }{\partial y^{ri}}\right\} ,$ $%
\frac{\partial }{\partial \overline{z}^{ri}}=\frac{1}{2}\left\{ \frac{%
\partial }{\partial x^{ri}}+{\mathbf{i}}\frac{\partial }{\partial y^{ri}}%
\right\} ,$ $dz^{ri}=dx^{ri}+{\mathbf{i}}dy^{ri},$ $d\overline{z}%
^{ri}=dx^{ri}-{\mathbf{i}}dy^{ri}.$ Note that $\{\frac{\partial }{\partial
z^{ri}},\frac{\partial }{\partial \overline{z}^{ri}}\}$ and $\{dz^{ri},d%
\overline{z}^{ri}\}$ are bases of the tangent space $T_{p}(^{k}M)$ and of
the cotangent space $T_{p}^{\ast }(^{k}M),$ respectively. The endomorphism $%
J_{k}$ is given by
\begin{equation}
J_{k}(\frac{\partial }{\partial z^{ri}})={\mathbf{i}}\frac{\partial }{%
\partial z^{ri}},\text{ }J_{k}(\frac{\partial }{\partial \overline{z}^{ri}}%
)=-{\mathbf{i}}\frac{\partial }{\partial \overline{z}^{ri}}.  \label{1.4}
\end{equation}%
If $J_{k}^{\ast }$ is an endomorphism of the cotangent space $T_{p}^{\ast
}(^{k}M)$ such that $J_{k}^{\ast 2}=-I$ $,$ then it holds
\begin{equation}
J_{k}^{\ast }(dz^{ri})={\mathbf{i}}dz^{ri},\text{ }J_{k}^{\ast }(d\overline{z%
}^{ri})=-{\mathbf{i}}d\overline{z}^{ri}.  \label{1.5}
\end{equation}

\subsection{\textbf{Higher Order Lifts of Complex Functions}}

In this subsection, we give definitions about higher order vertical and
complete lifts of complex functions defined on any complex manifold $M$ to $%
k $-th order extension $^{k}M$. The \textit{vertical lift} of a function $f$
defined on $M$ to $^{k}M$ is the function $f^{v^{k}}$ on $^{k}M$ given by
the equality:
\begin{equation}
f^{v^{k}}=f\circ \tau _{M}\circ \tau _{^{2}M}\circ ...\circ \tau _{^{k-1}M},
\label{1.6}
\end{equation}%
where $\tau _{^{k-1}M}:^{k}M\rightarrow ^{k-1}M$ is a canonical projection.
The \textit{complete lift} of function $f$ to $^{k}M$ is the function $%
f^{c^{k}}$ denoted by
\begin{equation}
f^{c^{k}}=\overset{.}{z}^{ri}(\frac{\partial f^{c^{k-1}}}{\partial z^{ri}}%
)^{v}+\overset{.}{\overline{z}}^{ri}(\frac{\partial f^{c^{k-1}}}{\partial
\overline{z}^{ri}})^{v}.  \label{1.7}
\end{equation}%
Using the induction method, the properties about vertical and complete lifts
of complex functions have been extended as follows:

$%
\begin{array}{ll}
i) & (f+g)^{v^{r}}=f^{v^{r}}+g^{v^{r}},(f.g)^{v^{r}}=f^{v^{r}}.g^{v^{r}} \\
ii) & (f+g)^{c^{r}}=f^{c^{r}}+g^{c^{r}},(f.g)^{c^{r}}=\sum_{j=0}^{r}\mathrm{C%
}_{j}^{r}f^{c^{r-j}v^{j}}.g^{c^{j}v^{r-j}}, \\
iii) & (\frac{\partial f}{\partial z^{0i}})^{v^{r}}=\frac{\partial f^{c^{r}}%
}{\partial z^{ri}},(\frac{\partial f}{\partial \overline{z}^{0i}})^{v^{r}}=%
\frac{\partial f^{c^{r}}}{\partial \overline{z}^{ri}}, \\
vi) & (\frac{\partial f}{\partial z^{0i}})^{c^{r}}=\frac{\partial f^{c^{r}}}{%
\partial z^{0i}},(\frac{\partial f}{\partial \overline{z}^{0i}})^{c^{r}}=%
\frac{\partial f^{c^{r}}}{\partial \overline{z}^{0i}},%
\end{array}%
$

where $f$ and $g$ are complex functions, and $\mathrm{C}_{j}^{r}$ is the
combination$.$

\subsection{\textbf{Higher Order Lifts of Complex Vector Fields}}

Here, the definitions and propositions about higher order vertical and
complete lifts of complex vector fields defined on any complex manifold $M$
to $k$-th order extension $^{k}M$ are presented. The \textit{vertical lift }%
of $X$ to $^{k}M$ is the complex vector field $X^{v^{k}}$ on $^{k}M$
formulated as below:
\begin{equation}
X^{v^{k}}(f^{c^{k}})=(Xf)^{v^{k}}.  \label{1.8}
\end{equation}

Now, we give the local expression of the vertical lift of $X$ to $^{k}M$.

\textbf{Proposition\thinspace 2.3.1} Let $M$ be any complex manifold and $%
^{k}M$ its $k$-th order extension$.$ Consider $X=Z^{0i}\frac{\partial }{%
\partial z^{0i}}+\overline{Z}^{0i}\frac{\partial }{\partial \overline{z}^{0i}%
}.$ Then the vertical lift of $X$ to $^{k}M$ is
\begin{equation}
X^{v^{k}}=(Z^{0i})^{v^{k}}\frac{\partial }{\partial z^{ki}}+(\overline{Z}%
^{0i})^{v^{k}}\frac{\partial }{\partial \overline{z}^{ki}}.  \label{1.9}
\end{equation}%
The \textit{complete lift} of $X$ to $^{k}M$ is the complex vector field $%
X^{c^{k}}$ such that
\begin{equation}
X^{c^{k}}(f^{c^{k}})=(Xf)^{c^{k}}.  \label{1.10}
\end{equation}%
The local expression of the complete lift of $X$ to $^{k}M$ is obtained as
follows.

\textbf{Proposition\thinspace 2.3.2 }Let $M$ be any complex manifold and $%
^{k}M$ its $k$-th order extension$.$ Let $X=Z^{0i}\frac{\partial }{\partial
z^{0i}}+\overline{Z}^{0i}\frac{\partial }{\partial \overline{z}^{0i}}$ .
Then the complete lift of $X$ to $^{k}M$ is
\begin{equation}
X^{c^{k}}=\mathrm{C}_{j}^{r}(Z^{0i})^{v^{k-r}c^{r}}\frac{\partial }{\partial
z^{ri}}+\mathrm{C}_{j}^{r}(\overline{Z}^{0i})^{v^{k-r}c^{r}}\frac{\partial }{%
\partial \overline{z}^{ri}}.  \label{1.11}
\end{equation}%
The extended properties about vertical and complete lifts of complex vector
fields by using the induction method are formulated as follows:

$%
\begin{array}{ll}
i) & (X+Y)^{v^{r}}=X^{v^{r}}+Y^{v^{r}},(X+Y)^{c^{r}}=X^{c^{r}}+Y^{c^{r}}, \\
ii) & (fX)^{v^{r}}=f^{v^{r}}X^{v^{r}},(fX)^{c^{r}}=\sum_{j=0}^{r}\mathrm{C}%
_{j}^{r}\,\,f^{c^{r-j}v^{j}}X^{c^{j}v^{r-j}}, \\
iii) &
\begin{array}{l}
X^{v^{k}}(f^{v^{k}})=0,\text{ }X^{c^{k}}(f^{c^{k}})=(Xf)^{c^{k}}, \\
X^{c^{k}}(f^{v^{k}})=X^{v^{k}}(f^{c^{k}})=(Xf)^{v^{k}}%
\end{array}
\\
iv) &
\begin{array}{l}
\left[ X^{v^{k}},Y^{v^{k}}\right] =0,\left[ X^{c^{k}},Y^{c^{k}}\right] =%
\left[ X,Y\right] ^{c^{k}}, \\
\left[ X^{v^{k}},Y^{c^{k}}\right] =\left[ X^{c^{k}},Y^{v^{k}}\right] =\left[
X,Y\right] ^{v^{k}}%
\end{array}
\\
v) & \Im _{0}^{1}(M)=Sp\left\{ \frac{\partial }{\partial z^{0i}},\frac{%
\partial }{\partial \overline{z}^{0i}}\right\} ,\Im
_{0}^{1}(^{k}M)=Sp\left\{ \frac{\partial }{\partial z^{ri}},\frac{\partial }{%
\partial \overline{z}^{ri}}\right\} , \\
&
\begin{array}{l}
(\frac{\partial }{\partial z^{0i}})^{c^{r}}=\frac{\partial }{\partial z^{0i}}%
,\text{ }(\frac{\partial }{\partial \overline{z}^{0i}})^{c^{r}}=\frac{%
\partial }{\partial \overline{z}^{0i}},\, \\
(\frac{\partial }{\partial z^{0i}})^{v^{r}}=\frac{\partial }{\partial z^{ri}}%
,\text{ }(\frac{\partial }{\partial \overline{z}^{0i}})^{v^{r}}=\frac{%
\partial }{\partial \overline{z}^{ri}}%
\end{array}%
\end{array}%
$

denoting by $M$ complex manifold$,$ by $X$,$Y$\ complex vector fields and by
$f$ function.

\subsection{\textbf{Higher Order Lifts of Complex 1- Forms}}

This subsection covers definitions and propositions about higher order
vertical and complete lifts of complex 1-forms defined on any complex
manifold $M$ to $k$-th order extension $^{k}M.$ The \textit{vertical lift}
of $\alpha $ to $^{k}M$ is the complex 1-form $\alpha ^{v^{k}}$ on $^{k}M$
defined by
\begin{equation}
\alpha ^{v^{k}}(X^{c^{k}})=(\alpha X)^{v^{k}}.  \label{1.12}
\end{equation}%
Now, we state a proposition on the vertical lift of $\alpha $ to $^{k}M$.

\textbf{Proposition\thinspace 2.4.1} Let $M$ be any complex manifold and $%
^{k}M$ its $k$-th order extension$.$ Set $\alpha =\alpha _{0i}dz^{0i}+%
\overline{\alpha }_{0i}d\overline{z}^{0i}$. Then the vertical lift of $%
\alpha $ to $^{k}M$ is
\begin{equation}
\alpha ^{v^{k}}=(\alpha _{0i})^{v^{k}}dz^{0i}+(\overline{\alpha }%
_{0i})^{v^{k}}d\overline{z}^{0i}.  \label{1.13}
\end{equation}%
The \textit{complete lift }of $\alpha $ to $^{k}M$ is the complex 1-form $%
\alpha ^{c^{k}}$ on $^{k}M$ defined by
\begin{equation}
\alpha ^{c^{k}}(X^{c^{k}})=(\alpha X)^{c^{k}}.  \label{1.14}
\end{equation}%
Now, we state a proposition on the complete lift of $\alpha $ to $^{k}M.$

\textbf{Proposition\thinspace 2.4.2} Let $M$ be any complex manifold and $%
^{k}M$ its $k$-th order extended complex manifold$.$ Put $\alpha =\alpha
_{0i}dz^{0i}+\overline{\alpha }_{0i}d\overline{z}^{0i}$. Then the complete
lift of $\alpha $ to $^{k}M$ is
\begin{equation}
\alpha ^{c^{k}}=(\alpha _{0i})^{c^{k-r}v^{r}}dz^{ri}+(\overline{\alpha }%
_{0i})^{c^{k-r}v^{r}}d\overline{z}^{ri}.  \label{1.15}
\end{equation}%
Using the induction method, the extended properties about vertical and
complete lifts of complex 1-forms are given as follows:

$%
\begin{array}{ll}
i) & (\alpha +\lambda )^{v^{r}}=\alpha ^{v^{r}}+\lambda ^{v^{r}},(\alpha
+\lambda )^{c^{r}}=\alpha ^{c^{r}}+\lambda ^{c^{r}}, \\
ii) & (f\alpha )^{v^{r}}=f^{v^{r}}\alpha ^{v^{r}},(f\alpha
)^{c^{r}}=\sum_{j=0}^{r}\mathrm{C}_{j}^{r}f^{c^{r-j}v^{j}}\alpha
^{c^{j}v^{r-j}}, \\
iii) & \Im _{1}^{0}(M)=Sp\left\{ dz^{0i},d\overline{z}^{0i}\right\} ,\Im
_{1}^{0}(^{k}M)=Sp\left\{ dz^{ri},d\overline{z}^{ri}\right\} , \\
&
\begin{array}{l}
(dz^{0i})^{c^{r}}=dz^{ri},(d\overline{z}^{0i})^{c^{r}}=d\overline{z}^{ri},\,
\\
(dz^{0i})^{v^{r}}=dz^{0i},(d\overline{z}^{0i})^{v^{r}}=d\overline{z}^{0i}.%
\end{array}%
\end{array}%
$

\subsection{\textbf{Higher Order Lifts of Almost Complex Structure}}

This subsection presents the definitions about higher order vertical and
complete lifts of an almost complex structure defined on any complex
manifold $M$ to the extended complex manifold $^{k}M$. The \textit{vertical
lifts} of $J_{0}$ and $J_{0}^{\ast }$ to $^{k}M$ are respectively the
structures $J_{0}^{v^{k}}$ and $J_{0}^{\ast v^{k}}$on $^{k}M$ defined by
\begin{equation}
J_{0}^{v^{k}}(X^{c^{k}})=(J_{0}X)^{v^{k}},J_{0}^{\ast v^{k}}(\alpha
^{c^{k}})=(J_{0}^{\ast }\alpha )^{v^{k}}.  \label{1.16}
\end{equation}%
The \textit{complete lifts }of $J_{0}$ and $J_{0}^{\ast }$ to $^{k}M$ are
respectively the structures $J_{0}^{c^{k}}$ and $J_{0}^{\ast c^{k}}$on $%
^{k}M $ denoted by
\begin{equation}
J_{0}^{c^{k}}(X^{c^{k}})=(J_{0}X)^{c^{k}},J_{0}^{\ast c^{k}}(\alpha
^{c^{k}})=(J_{0}^{\ast }\alpha )^{c^{k}}.  \label{1.17}
\end{equation}

\subsection{\textbf{Higher Order Lifts of Hermitian Metric}}

This subsection gives the definitions about higher order vertical and
complete lifts of a Hermitian metric on any complex manifold $M$ to the
extension space $^{k}M$. In order to define higher order vertical and
complete lifts of a Hermitian metric on $M$, firstly the definition of
higher order vertical and complete lifts of complex tensor fields of type
(0,2) is given. The \textit{vertical lift} of $G$ to $^{k}M$ is the tensor
field of type (0,2) $G^{v^{k}}$on $^{k}M$ formulated as
\begin{equation}
G^{v^{k}}(X^{c^{k}},Y^{c^{k}})=(G(X,Y))^{v^{k}}.  \label{1.18}
\end{equation}%
Denote by $g$ a Hermitian metric and by $J_{0}$ an almost complex structure
on any complex manifold $M.$ Since $g$ is a tensor field of type (0,2)$,$ we
have the equality
\begin{equation}
g^{v^{k}}(X^{c^{k}},Y^{c^{k}})=(g(X,Y))^{v^{k}}=g^{v^{k}}(J_{0}^{c^{k}}X^{c^{k}},J_{0}^{c^{k}}Y^{c^{k}}),
\label{1.19}
\end{equation}%
for any complex vector fields $X,Y$ on $M.$ Hence the \textit{vertical lift}
of $g$ to $^{k}M$ \ is a Hermitian metric $g^{v^{k}}$on $^{k}M.$ The
extended complex manifold $^{k}M$ with Hermitian metric $g^{v^{k}}$ is
called the \textit{vertical lift} of order $k$ of the Hermitian manifold $%
(M,J_{0},g).$ The \textit{complete lift }of $G$ to $^{k}M$ is the tensor
field of type (0,2) $G^{c^{k}}$ on $^{k}M$ given by equality
\begin{equation}
G^{c^{k}}(X^{c^{k}},Y^{c^{k}})=(G(X,Y))^{c^{k}}.  \label{1.20}
\end{equation}%
Also we find\ the equality
\begin{equation}
g^{c^{k}}(X^{c^{k}},Y^{c^{k}})=(g(X,Y))^{c^{k}}=g^{c^{k}}(J_{0}^{c^{k}}X^{c^{k}},J_{0}^{c^{k}}Y^{c^{k}}),
\label{1.21}
\end{equation}%
for any complex vector fields $X,Y$ on $M.$ Hence the \textit{complete}
\textit{lift} of $g$ to $^{k}M$ is $g^{c^{k}}$on $^{k}M.$ The extended
complex manifold $^{k}M$ with Hermitian metric $g^{c^{k}}$ is called \textit{%
complete lift} of order $k$ of Hermitian manifold$.$

\subsection{\textbf{Higher Order Lifts of Kaehlerian Form}}

In this part, the higher order vertical and complete lifts of a Kaehlerian
form on any complex manifold $M$ to its extension $^{k}M$ are defined. Given
a Hermitian manifold $(M,J_{0},g)$ . Since the Kaehlerian form $\Phi $ is
tensor field of type (0,2)$,$ we obtain the equality
\begin{equation}
\Phi ^{v^{k}}(X^{c^{k}},Y^{c^{k}})=(\Phi
(X,Y))^{v^{k}}=g^{v^{k}}(X^{c^{k}},J_{0}^{c^{k}}Y^{c^{k}}),  \label{1.22}
\end{equation}%
for any complex vector fields $X,Y$ on $M.$ Hence the \textit{vertical lift}
of $\Phi $ to $^{k}M$ is a Kaehlerian form $\Phi ^{v^{k}}$on $^{k}M.$ Also
we have equality
\begin{equation}
\Phi ^{c^{k}}(X^{c^{k}},Y^{c^{k}})=(\Phi
(X,Y))^{c^{k}}=g^{c^{k}}(X^{c^{k}},J_{0}^{c^{k}}Y^{c^{k}}),  \label{1.23}
\end{equation}%
for any complex vector fields $X,Y$ on $M.$ Thus the \textit{complete lift}
of $\Phi $ to extended complex manifold $^{k}M$ is a Kaehlerian form $\Phi
^{c^{k}}$on $^{k}M.$

\subsection{\textbf{Higher Order Lifts of Kaehlerian Metric}}

Higher order vertical and complete lifts of a Kaehlerian form associated
with any Hermitian manifold $M$ to its extension $^{k}M$ is introduced in
this subsection. Also we give definitions about higher order vertical and
complete lifts of a Kaehlerian metric defined on $M$ to $^{k}M$.

One needs to specify in the statement that\textbf{\ }$(M,J_{0},g)$ is a
Kaehlerian manifold, since only in this case $d\Phi =0.$ In fact some
authors call Kaehlerian form the 2-form $\Phi $ associated with any
Hermitian manifold $M.$ Let $\Phi ^{v^{k}}$ (resp. $\Phi ^{c^{k}}$) its $k$%
-th order vertical lift ( resp. complete lift). Then we have
\begin{equation}
d\Phi ^{v^{k}}=0,\text{ }d\Phi ^{c^{k}}=0.  \label{1.24}
\end{equation}%
Since the Kaehlerian form $\Phi ^{v^{k}}($resp. $\Phi ^{c^{k}})$ on $^{k}M$
is closed, the Hermitian metric $g^{v^{k}}$ (resp. $g^{c^{k}}$) on $^{k}M$
is said to be the \textit{vertical lift }(resp.\textit{\ complete lift}) of
Kaehlerian metric $g$ to $^{k}M$. The extended Hermitian manifold $^{k}M$
with Kaehlerian metric $g^{v^{k}}$ (resp. $g^{c^{k}}$) is called \textit{%
vertical lift }(resp. \textit{complete lift}) of order $k$ of Kaehlerian
manifold.

\section{\textbf{Product Manifolds and Lifted Structures}}

From this point onwards the definitions and structures given in \cite%
{bowman, civelek,tekkoyun1} can be extended as follows.

\textbf{Definition 3.1:} Let $^{0}N=\mathbf{R}\times ^{0}M,$\ $^{1}N=\mathbf{%
R}\times ^{1}M$\ and $^{2}N$\ $=\mathbf{R}\times ^{2}M$\ be manifolds.
Consider a sequence given by
\begin{equation}
^{0}\mathbf{N}\overset{\tau _{^{0}N}}{\longleftarrow }\text{ }^{1}\mathbf{N}%
\overset{\tau _{^{1}N}}{\longleftarrow }\text{ \ }^{2}\mathbf{N~}\text{\ }
\label{2.0}
\end{equation}%
where $\tau _{^{0}N}\ $and $\tau _{^{1}N}$\ are smooth maps. If the kernel
of the map $\tau _{^{0}N}$\ is equal to the image set of the map $\tau
_{^{1}N},$\ then the sequence (\ref{2.0}) is said to be a \textit{short
exact sequence}.

\textbf{Definition 3.2: }If $N$\ is a manifold, then a sequence of $N$\ is a
sequence $S$\ of manifolds and maps determined by sequence
\begin{equation}
^{0}N\overset{\tau _{^{0}N}}{\longleftarrow }\text{ }^{1}N\overset{\tau
_{^{1}N}}{\longleftarrow }\text{ \ }^{2}N~\overset{\tau _{^{2}N}}{%
\longleftarrow }...\text{\ }  \label{2.01}
\end{equation}%
where $^{0}N=N$ and each short sequence of the sequence (\ref{2.01}) is an
\textit{exact sequence} of $N$. If $S$\ has a last term $^{m}N$\ \ then we
say that $S$\ has length $m$, otherwise we say that it has infinite length
or length $\infty .$ Denote by $T(\tau _{^{k}N})$\ the differential (tangent
functor) of $\tau _{^{k}N},$\ by $T(^{k}N)$\ the tangent bundle of $^{k}N$\
and by $\tau _{^{k}N}$\ the natural projection of $T(^{k}N)$\ to $^{k}N$\ .

\textbf{Definition 3.3:} If $N$\ is a manifold, a sequence $S$\ of length $%
m\leq \infty $\ of $N$\ yields the following properties: \newline
i) For each integer $1\leq k\leq m$\ there exists on imbedding $^{k-1}I$:$%
^{k}N\rightarrow T(^{k-1}N)$\ with $^{0}I$\ onto such that $\tau
_{^{k-1}N}=\theta _{^{k}N}\circ ^{k-1}I;$ \newline
ii) For each integer $1\leq k\leq m$\ \ the diagram%
\begin{equation*}
\begin{array}{ccc}
T(^{k}N) & \overset{T(\theta _{^{k-1}N})}{\longrightarrow } & T(^{k-1}N) \\
\theta _{^{k}N}\downarrow &  & \downarrow I \\
^{k}N & \overset{^{k-1}I}{\longrightarrow } & T(^{k-1}N)%
\end{array}%
\end{equation*}%
commutes on $^{k}I(^{k+1}M)$\ exactly. In this case $S$\ is said to be an
extended sequence of length $m$\ of $N$, and $^{k}N=\mathbf{R}\times ^{k}M$\
\ is called a $k$-\textit{th order extension of product manifold} $N=\mathbf{%
R}\times M$\ $($or $^{0}N=\mathbf{R}\times ^{0}M)$\ of dimension $2m+1,$\
where $^{k}M$\ is an extended complex manifold.

Let $(t,z^{ri},\,\overline{z}^{ri})$ be a coordinate system on a
neighborhood $^{k}V$ of any point $p$ of $^{k}N$. Therefore, we respectively
define by$\left\{ \frac{\partial }{\partial t},\frac{\partial }{\partial
z^{ri}},\frac{\partial }{\partial \overline{z}^{ri}}\right\} $ and $\left\{
dt,dz^{ri},d\overline{z}^{ri}\right\} $ the natural bases over coordinate
system of tangent space $T_{p}(^{k}N)$ and cotangent space $T_{p}^{\ast
}(^{k}N)$ of $^{k}N.$

Let $f$ be a function defined on $N$ and $(t,z^{0i},\overline{z}^{0i})$ be
coordinates of $N.$ So, the $1$- form defined by
\begin{equation}
df=\frac{\partial f}{\partial t}dt+\frac{\partial f}{\partial z^{0i}}dz^{0i}+%
\frac{\partial f}{\partial \overline{z}^{0i}}d\overline{z}^{0i}  \label{2.1}
\end{equation}%
is the differential of $f$. Denote by $\chi (N)$ the set of vector fields
and by $\chi ^{\ast }(N)$ the set of dual vector fields on $N.$ In this
case, elements $Z$ and $\omega $ of $\chi (N)$ and $\chi ^{\ast }(N)$ are
determined by
\begin{equation}
Z=\frac{\partial }{\partial t}+Z^{0i}\frac{\partial }{\partial z^{0i}}+%
\overline{Z}^{0i}\frac{\partial }{\partial \overline{z}^{0i}}  \label{2.2}
\end{equation}%
and
\begin{equation}
\omega =dt+\omega _{0i}dz^{0i}+\overline{\omega }_{0i}d\overline{z}^{0i},
\label{2.3}
\end{equation}%
respectively, where $Z^{0i},\omega _{0i},\overline{Z}^{0i},\overline{\omega }%
_{0i}\in \mathcal{F}(M)$.

\subsection{\textbf{Higher Order Lifts of Functions}}

In this subsection, extensions of definitions and properties about the
higher order vertical, complete, complete-vertical and horizontal lifts of
functions on product manifold $N$ of dimension $2m+1$ to $^{k}N$ are
obtained.

Let $^{k-1}N$ be the $(k-1)$-th order extension of $N$ and $\tau
_{^{k-1}N}:^{k}N\rightarrow ^{k-1}N$ the natural projection. Consider the
linear isomorphism as follows:
\begin{equation}
\begin{array}{llll}
v: & \Im _{0}^{0}(^{k-1}N) & \rightarrow  & \Im _{0}^{0}(^{k}N) \\
& \widetilde{f} & \rightarrow  & v(\widetilde{f})=\widetilde{f}^{v}.%
\end{array}
\label{3.1}
\end{equation}%
Thus, the \textit{vertical lift} of function $\widetilde{f}$ to $^{k}N$ is
the function $\widetilde{f}^{v}=\widetilde{f}\circ \tau _{^{k-1}N}$ .

Let $f^{v^{k-1}}$ be vertical lift of a function $f$ on $N$ to $^{k-1}N.$ In
(\ref{3.1}), if $\widetilde{f}=f^{v^{k-1}},$ then the \textit{vertical lift}
of function $f$ to $^{k}N$ is the function $f^{v^{k}}$ defined by the
equality%
\begin{equation}
f^{v^{k}}=f\circ \tau _{N}\circ \tau _{^{2}N}\circ ...\circ \tau _{^{k-1}N}.
\label{3.2}
\end{equation}%
Let $\widetilde{f}$ be a function on $^{k-1}N$ and $(t,z^{ri},\overline{z}%
^{ri}),0\leq r\leq k-1$ be extended coordinates of $^{k-1}N.$ Therefore, the
1- form defined by the equality
\begin{equation}
d\widetilde{f}=\frac{\partial \widetilde{f}}{\partial t}dt+\frac{\partial
\widetilde{f}}{\partial z^{ri}}dz^{ri}+\frac{\partial \widetilde{f}}{%
\partial \overline{z}^{ri}}d\overline{z}^{ri}  \label{3.3}
\end{equation}%
is a differential of $\widetilde{f}.$

Suppose that let $\iota _{k}$ be the linear isomorphism such that
\begin{equation}
\iota _{k}:\Im _{1}^{0}(^{k-1}N)\rightarrow \Im _{0}^{0}(^{k}N)  \label{3.4}
\end{equation}%
\begin{equation*}
\iota _{k}(dt)=t,\iota _{k}(dz^{ri})=\overset{.}{z}^{ri},\,\iota _{k}(d%
\overline{z}^{ri})=\overset{.}{\overline{z}}^{ri},
\end{equation*}%
where $Sp\left\{ dt,dz^{ri},d\overline{z}^{ri}:0\leq r\leq k-1\right\} =\Im
_{1}^{0}(^{k-1}N).$ Taking account of (\ref{3.3}), we say that the \textit{%
complete lift} of the function $\widetilde{f}$ to $^{k}N$ is the function $%
\widetilde{f}^{c}$ defined by
\begin{equation}
\widetilde{f}^{c}=\iota _{k}(d\widetilde{f})=t(\frac{\partial \widetilde{f}}{%
\partial t})^{v}+\overset{.}{z}^{ri}(\frac{\partial \widetilde{f}}{\partial
z^{ri}})^{v}+\overset{.}{\overline{z}}^{ri}(\frac{\partial \widetilde{f}}{%
\partial \overline{z}^{ri}})^{v}.  \label{3.5}
\end{equation}%
Let $f^{c^{k-1}}$ be complete lift of a function $f$ on $N$ to $^{k-1}N.$ In
(\ref{3.5}), if $\widetilde{f}=f^{c^{k-1}},$ the \textit{complete lift} of
function $f$ to extended manifold $^{k}N$ is the function $f^{c^{k}}$
defined by equality
\begin{equation}
f^{c^{k}}=t(\frac{\partial f^{c^{k-1}}}{\partial t})^{v}+\overset{.}{z}^{ri}(%
\frac{\partial f^{c^{k-1}}}{\partial z^{ri}})^{v}+\overset{.}{\overline{z}}%
^{ri}(\frac{\partial f^{c^{k-1}}}{\partial \overline{z}^{ri}})^{v}.
\label{3.6}
\end{equation}%
Let $f^{c^{r}}$ be $r$ -th order complete lift of a function $f$ on $N$ to $%
^{r}N.$ Then if it is taken $s$-th order vertical lift of function $%
f^{c^{r}} $ on $^{r}N$ to $^{k}N,$ by \textit{complete}-\textit{vertical
lift of order }$(r,s)$\textit{\ }of $f$ on $N$ to $^{k}N$ we call the
function $f^{c^{r}v^{s}}$determined by
\begin{equation}
(f^{c^{r}})^{v^{s}}=f^{c^{r}v^{s}}=f^{c^{r}}\circ \tau _{^{r}N}\circ
...\circ \tau _{^{r+s-1}N},  \label{3.7}
\end{equation}%
where $0\leq r,s\leq k$ and $r+s=k.$

There exists commutative property taking into complete-vertical lift of
functions. i.e., it means that complete-vertical lifts of order $(r,s)$ and
complete-vertical lifts of order $(s,r)$ of functions on $N$ to its
extension $^{k}N$ are equal. The \textit{horizontal\thinspace \thinspace lift%
} of\thinspace \thinspace $f$ on $N$ \thinspace to $^{k}N$ is the function $%
f^{H^{k}}$ on $^{k}N$ given by
\begin{equation}
f^{H^{k}}=f^{c^{k}}-\gamma (\nabla f^{c^{k-1}}),\,\,\,\,\,(\gamma (\nabla
f^{c^{k-1}})=\nabla _{\gamma }f^{c^{k-1}}),  \label{3.8}
\end{equation}%
where $\nabla $ is an affine linear connection on $^{k-1}N$ with local
components $\Gamma _{rj}^{ri},1\leq i,j\leq m,$ $\nabla f^{c^{k-1}}$ is
gradient of $f^{c^{k-1}}$ and $\gamma $ is an operator given by
\begin{equation*}
\gamma :\Im _{s}^{r}(^{k-1}N)\rightarrow \Im _{s-1}^{r}(^{k}N).
\end{equation*}%
Thus, one has $f^{H^{k}}=0,$ since
\begin{equation*}
\nabla _{\gamma }f^{c^{k-1}}=(\frac{\partial f^{c^{k-1}}}{\partial t})^{v}+%
\overset{.}{z}^{ri}(\frac{\partial f^{c^{k-1}}}{\partial z^{ri}})^{v}+%
\overset{.}{\overline{z}}^{ri}(\frac{\partial f^{c^{k-1}}}{\partial
\overline{z}^{ri}})^{v},
\end{equation*}%
where dots mean derivation with respect to time. The generic properties of
the higher order vertical, complete and horizontal lifts of functions on $N$
are

$%
\begin{array}{ll}
i) &
\begin{array}{l}
(f+g)^{v^{r}}=f^{v^{r}}+g^{v^{r}},(f+g)^{H^{k}}=0, \\
(f+g)^{c^{r}}=f^{c^{r}}+g^{c^{r}}%
\end{array}
\\
ii) &
\begin{array}{l}
(f.g)^{v^{r}}=f^{v^{r}}.g^{v^{r}},(f.g)^{H^{k}}=0, \\
(f.g)^{c^{r}}=\sum_{j=0}^{r}\mathrm{C}%
_{j}^{r}f^{c^{r-j}v^{j}}.g^{c^{j}v^{r-j}}%
\end{array}
\\
iii) &
\begin{array}{l}
(\frac{\partial f}{\partial z^{0i}})^{v^{r}}=\frac{\partial f^{c^{r}}}{%
\partial z^{ri}},(\frac{\partial f}{\partial \overline{z}^{0i}})^{v^{r}}=%
\frac{\partial f^{c^{r}}}{\partial \overline{z}^{ri}}, \\
(\frac{\partial f}{\partial t})^{v^{r}}=\frac{\partial f^{c^{r}}}{\partial t}%
\end{array}
\\
iv) &
\begin{array}{l}
(\frac{\partial f}{\partial z^{0i}})^{c^{r}}=\frac{\partial f^{c^{r}}}{%
\partial z^{0i}},(\frac{\partial f}{\partial \overline{z}^{0i}})^{c^{r}}=%
\frac{\partial f^{c^{r}}}{\partial \overline{z}^{0i}}, \\
(\frac{\partial f}{\partial t})^{c^{r}}=\frac{\partial f^{c^{r}}}{\partial t}%
\end{array}%
\end{array}%
$

for all $f,g\in \mathcal{F}($ $N)$.

\subsection{\textbf{Higher Order Lifts of Vector Fields}}

Extensions of definitions and propositions about higher order vertical,
complete, complete-vertical and horizontal lifts of vector fields on product
manifold $N$ of dimension $2m+1$ to $^{k}N$ are derived in this subsection.

Let $^{k-1}N$ be the $(k-1)$-th order extension of $N.$ Denote by $%
\widetilde{Z}$ a vector field on $^{k-1}N.$ Then the \textit{vertical lift}
of the vector field $\widetilde{Z}$ to $^{k}N$ is the vector field $%
\widetilde{Z}^{v}$ on $^{k}N$ defined by:
\begin{equation}
\widetilde{Z}^{v}(\widetilde{f}^{c})=(\widetilde{Z}\widetilde{f})^{v}.
\label{4.1}
\end{equation}%
We denote by $f^{c^{k-1}}$ and $Z^{v^{k-1}}$complete lift of a function $f$
and the vertical lift of a vector field $Z$ on $N$ to $^{k-1}N,$ respectively%
$.$ In (\ref{4.1}), if $\widetilde{f}=f^{c^{k-1}}$ and $\widetilde{Z}%
=Z^{v^{k-1}},$ then the \textit{vertical lift }of vector field $Z$ on $N$ to
$^{k}N$ is the vector field $Z^{v^{k}}$ on $^{k}N$ given by%
\begin{equation}
Z^{v^{k}}(f^{c^{k}})=(Zf)^{v^{k}}.  \label{4.2}
\end{equation}%
\textbf{Proposition\thinspace 3.2.1} Let $^{k}N$ be $k$-th extension of $N.$
Assume that the vector field $Z$ defined on $N$ is given by (\ref{2.2})$.$
Then vertical lift of order $k$ of $Z$ to $^{k}N$ is
\begin{equation}
Z^{v^{k}}=\frac{\partial }{\partial t}+(Z^{0i})^{v^{k}}\frac{\partial }{%
\partial z^{ki}}+(\overline{Z}^{0i})^{v^{k}}\frac{\partial }{\partial
\overline{z}^{ki}}.  \label{4.3}
\end{equation}%
\textbf{Proof: }Considering a coordinate system $(t,z^{ri},\,\overline{z}%
^{ri})$ on a neighborhood $^{k}V$ of any point $p$ of $^{k}N$, we put $%
Z^{v^{k}}=\frac{\partial }{\partial t}+Z^{ri}\frac{\partial }{\partial z^{ri}%
}+\overline{Z}^{ri}\frac{\partial }{\partial \overline{z}^{ri}}.$ Let $%
f^{c^{k}}$ be $k$-th order complete lift of function $f$ to $^{k}N.$ Then
from vertical lift properties we can obtain equations
\begin{equation}
Z^{v^{k}}(f^{c^{k}})=\frac{\partial f^{c^{k}}}{\partial t}+Z^{ri}\frac{%
\partial f^{c^{k}}}{\partial z^{ri}}+\overline{Z}^{ri}\frac{\partial
f^{c^{k}}}{\partial \overline{z}^{ri}}  \label{4.3.1}
\end{equation}%
and
\begin{eqnarray}
(Zf)^{v^{k}} &=&(\frac{\partial f}{\partial t}+Z^{0i}\frac{\partial f}{%
\partial z^{0i}}+\overline{Z}^{0i}\frac{\partial f}{\partial \overline{z}%
^{0i}})^{v^{k}}  \label{4.3.2} \\
&=&\frac{\partial f^{c^{k}}}{\partial t}+(Z^{0i})^{v^{k}}\frac{\partial
f^{c^{k}}}{\partial z^{ki}}+(\overline{Z}^{0i})^{v^{k}}\frac{\partial
f^{c^{k}}}{\partial \overline{z}^{ki}}.  \notag
\end{eqnarray}%
By (\ref{4.2}), (\ref{4.3.1}), (\ref{4.3.2}) one obtains%
\begin{equation*}
Z^{ri}=0,\overline{Z}^{ri}=0,Z^{ki}=(Z^{0i})^{v^{k}},\overline{Z}^{ki}=(%
\overline{Z}^{0i})^{v^{k}},0\leq r\leq k-1.
\end{equation*}%
Hence, the proof finishes.$\Box $

Let $^{k-1}N$ be the $(k-1)$-th order extended manifold of $N$ and $%
\widetilde{Z}$ be a vector field on $^{k-1}N.$ Then the \textit{complete lift%
} of $\widetilde{Z}$ on $^{k-1}N$ to $^{k}N$ is defined by:
\begin{equation}
\widetilde{Z}^{c}(\widetilde{f}^{c})=(\widetilde{Z}\widetilde{f})^{c},
\label{4.4}
\end{equation}%
for any function $\widetilde{f}$ on $^{k-1}N.$  Let $f^{c^{k-1}}$ and $%
Z^{c^{k-1}}$be respectively complete lifts of a function $f$ and a vector
field $Z$ defined on $N$ to $^{k-1}N.$ In (\ref{4.4}), if $\widetilde{f}%
=f^{c^{k-1}}$ and $\widetilde{Z}=Z^{c^{k-1}},$ then the \textit{complete
lift }of vector field $Z$ on $N$ to extended manifold $^{k}N$ is the vector
field $Z^{c^{k}}$ on $^{k}N$ given by
\begin{equation}
Z^{c^{k}}(f^{c^{k}})=(Zf)^{c^{k}},  \label{4.5}
\end{equation}%
for any function $\widetilde{f}$ on $^{k-1}N$ .

Similar to the proof of \textbf{Proposition 3.2.1,} one may prove the
following:

\textbf{Proposition\thinspace \thinspace 3.2.2 }Let $^{k}N$ be extension of
order $k$ of $N.$ Assume that the vector field $Z$ on $N$ is given by (\ref%
{2.2})$.$ Then $k$-th order complete lift of $Z$ to $^{k}N$ is
\begin{equation}
Z^{c^{k}}=\frac{\partial }{\partial t}+\mathrm{C}%
_{r}^{k}(Z^{0i})^{v^{k-r}c^{r}}\frac{\partial }{\partial z^{ri}}+\mathrm{C}%
_{r}^{k}(\overline{Z}^{0i})^{v^{k-r}c^{r}}\frac{\partial }{\partial
\overline{z}^{ri}}.  \label{4.6}
\end{equation}%
Let $Z$ be a vector field on manifold $N.$ Then the \textit{complete}-%
\textit{vertical lift of order }$(r,s)$ of $Z$ to $^{k}N$ is the vector
field $Z^{c^{r}v^{s}}\in \chi (^{k}N)$ determined by equality
\begin{equation}
Z^{c^{r}v^{s}}(f^{c^{k}})=(Zf)^{c^{r}v^{s}},  \label{4.7}
\end{equation}%
where $0\leq r,s\leq k$ and $r+s=k.$

\textbf{Proposition\thinspace 3.2.3} Let $N$ be any product manifold of
dimension $2m+1$ and $^{k}N$ its $k$-th order extended manifold$.$ Consider
the vector field $Z$ on $N$ given by (\ref{2.2})$.$ Then the
complete-vertical lift of order $(r,s)$ of $Z$ to $^{k}N$ is
\begin{equation*}
Z^{c^{r}v^{s}}:\left(
\begin{array}{l}
\,\,\,\,\,\,\,\,\,\,\,\,\,\,\,\,1 \\
\mathrm{C}_{k-1}^{r}(Z^{0i})^{v^{s+k-l}c^{l-s}} \\
\mathrm{C}_{k-1}^{r}(\overline{Z}^{0i})^{v^{s+k-l}c^{l-s}}%
\end{array}%
\right) ,0\leq l\leq k.
\end{equation*}%
\textbf{Proof: }Considering a coordinate system $(t,z^{li},\,\overline{z}%
^{li})$ on a neighborhood $^{k}V$ of any point $p$ of $^{k}N$, we put $%
Z^{c^{r}v^{s}}=\frac{\partial }{\partial t}+Z^{li}\frac{\partial }{\partial
z^{li}}+\overline{Z}^{li}\frac{\partial }{\partial \overline{z}^{li}}$. Let $%
f^{c^{k}}$ be $k$-th order complete lift of function $f$ to the extended
manifold $^{k}N.$ Then from complete and vertical lift properties we
calculate
\begin{equation}
Z^{c^{r}v^{s}}(f^{c^{k}})=\frac{\partial f^{c^{k}}}{\partial t}+Z^{hi}\frac{%
\partial f^{c^{k}}}{\partial z^{li}}+\overline{Z}^{hi}\frac{\partial
f^{c^{k}}}{\partial \overline{z}^{li}}  \label{4.7.1}
\end{equation}%
and

\begin{eqnarray}
(Zf)^{c^{r}v^{s}} &=&(\frac{\partial f}{\partial t}+Z^{0i}\frac{\partial f}{%
\partial z^{0i}}+\overline{Z}^{0i}\frac{\partial f}{\partial \overline{z}%
^{0i}})^{c^{r}v^{s}}  \label{4.7.2} \\
&=&\frac{\partial f^{c^{k}}}{\partial t}+\mathrm{C}%
_{h}^{r}(Z^{0i})^{v^{s+h}c^{r-h}}\frac{\partial f^{c^{k}}}{\partial z^{k-hi}}
\notag \\
&&+\mathrm{C}_{h}^{r}(\overline{Z}^{0i})^{v^{s+h}c^{r-h}}\frac{\partial
f^{c^{k}}}{\partial \overline{z}^{k-hi}}.  \notag
\end{eqnarray}%
According to (\ref{4.7}), by (\ref{4.7.1}) and (\ref{4.7.2}) and also using $%
l=k-h$ from the following equalities
\begin{equation*}
\frac{\partial f^{c^{k}}}{\partial z^{li}}=\frac{\partial f^{c^{k}}}{%
\partial z^{k-hi}}\mathit{\ }\text{and }\frac{\partial f^{c^{k}}}{\partial
\overline{z}^{li}}=\frac{\partial f^{c^{k}}}{\partial \overline{z}^{k-hi}}
\end{equation*}%
we have
\begin{equation*}
Z^{li}=\mathrm{C}_{k-1}^{r}(Z^{0i})^{v^{s+k-l}c^{l-s}},\overline{Z}^{ri}=%
\mathrm{C}_{k-1}^{r}(\overline{Z}^{0i})^{v^{s+k-l}c^{l-s}},0\leq l\leq k
\end{equation*}%
Hence, the proof is completed.$\Box $

There exists commutative property considering complete-vertical lift of
vector fields. i.e., it means that complete-vertical lifts of order $(r,s)$
and complete-vertical lifts of order $(s,r)$ of vector fields on $N$ to its
extension $^{k}N$ are equal. The complete-vertical lifts of order $(r,s)$ of
vector fields on a manifold $N$ obey the following property
\begin{equation*}
(fZ)^{c^{r}v^{s}}=\mathrm{C}_{k-1}^{r}f^{v^{s+h}c^{r-h}}Z^{c^{h}v^{k-h}},\,0%
\leq r,s\leq k,\,(r+s=k).
\end{equation*}%
The \textit{horizontal lift} of a vector field\textbf{\ }$Z$ on $N$ to $%
^{k}N $ is the vector field $Z^{H^{k}}$ on $^{k}N$ given by
\begin{equation}
Z^{H^{k}}f^{v^{k}}=(Zf)^{v^{k}}.  \label{4.8}
\end{equation}%
Obviously, we have
\begin{equation*}
Z^{H^{k}}=\frac{\partial }{\partial t}+Z^{ri}D_{ri}+\overline{Z}^{ri}%
\overline{D}_{ri}
\end{equation*}%
such that $D_{ri}=\frac{\partial }{\partial z^{ri}}-\Gamma _{rj}^{ri}\frac{%
\partial }{\partial z^{r+1i}}$ and $\overline{D}_{ri}=\frac{\partial }{%
\partial \overline{z}^{ri}}-\overline{\Gamma }_{rj}^{ri}\frac{\partial }{%
\partial \overline{z}^{r+1i}},1\leq i,j\leq m.$ By an \textit{extended}
\textit{frame }adapted to a connection $\nabla $ on $^{k}N,$ we mean the set
of local vector fields $\left\{ \frac{\partial }{\partial t},D_{ri},%
\overline{D}_{ri},V_{ri}=\frac{\partial }{\partial z^{r+1i}},\overline{V}%
_{ri}=\frac{\partial }{\partial \overline{z}^{r+1i}}\right\} .$ The higher
order vertical, complete and horizontal lifts of vector fields on $N$ obey
the general properties

$%
\begin{array}{ll}
i) &
\begin{array}{l}
(Z+W)^{v^{r}}=Z^{v^{r}}+W^{v^{r}},(Z+W)^{c^{r}}=Z^{c^{r}}+W^{c^{r}}, \\
(Z+W)^{H^{k}}=Z^{H^{k}}+W^{H^{k}}%
\end{array}
\\
ii) & (fZ)^{v^{r}}=f^{v^{r}}Z^{v^{r}},(fZ)^{c^{r}}=\sum_{j=0}^{r}\mathrm{C}%
_{j}^{r}f^{c^{r-j}v^{j}}Z^{c^{j}v^{r-j}}, \\
iii) &
\begin{array}{l}
Z^{v^{k}}(f^{v^{k}})=0,Z^{c^{k}}(f^{v^{k}})=Z^{v^{k}}(f^{c^{k}})=(Zf)^{v^{k}},
\\
Z^{c^{k}}(f^{c^{k}})=(Zf)^{c^{k}},Z^{H^{k}}(f^{v^{k}})=(Zf)^{v^{k}}%
\end{array}
\\
iv) & \chi (V)=Sp\left\{ \frac{\partial }{\partial t},\frac{\partial }{%
\partial z^{0i}},\frac{\partial }{\partial \overline{z}^{0i}}\right\} ,\text{
}\chi (^{k}V)=Sp\left\{ \frac{\partial }{\partial t},\frac{\partial }{%
\partial z^{ri}},\frac{\partial }{\partial \overline{z}^{ri}}\right\} , \\
&
\begin{array}{l}
(\frac{\partial }{\partial z^{0i}})^{c^{r}}=\frac{\partial }{\partial z^{0i}}%
,(\frac{\partial }{\partial \overline{z}^{0i}})^{c^{r}}=\frac{\partial }{%
\partial \overline{z}^{0i}},(\frac{\partial }{\partial t})^{c^{r}}=\frac{%
\partial }{\partial t}, \\
(\frac{\partial }{\partial z^{0i}})^{v^{r}}=\frac{\partial }{\partial z^{ri}}%
,(\frac{\partial }{\partial \overline{z}^{0i}})^{v^{r}}=\frac{\partial }{%
\partial \overline{z}^{ri}}, \\
(\frac{\partial }{\partial t})^{v^{r}}=\frac{\partial }{\partial t},(\frac{%
\partial }{\partial t})^{H^{k}}=\frac{\partial }{\partial t}, \\
(\frac{\partial }{\partial z^{0i}})^{H^{k}}=D_{ri},(\frac{\partial }{%
\partial \overline{z}^{0i}})^{H^{k}}=\overline{D}_{ri},%
\end{array}%
\end{array}%
$

for all $Z,W\in \chi (N)$ and $f\in \mathcal{F}(N)$.

\subsection{\textbf{Higher Order Lifts of 1-Forms}}

In this subsection, we introduce definitions and propositions about higher
order vertical, complete, complete-vertical and horizontal lifts of a 1-form
defined on a product manifold $N$ of dimension $2m+1$ to $^{k}N$.

Let $^{k-1}N$ be $(k-1)$-th order extension of $N.$ Given a 1-form $%
\widetilde{\omega }$ on $^{k-1}N$ \ the \textit{vertical lift} of $%
\widetilde{\omega }$ to $^{k}N$ is the 1-form $\widetilde{\omega }^{v}$ on $%
^{k}N$ defined by:
\begin{equation}
\widetilde{\omega }^{v}(\widetilde{Z}^{c})=(\widetilde{\omega }\widetilde{Z}%
)^{v},  \label{5.1}
\end{equation}%
for any vector field $\widetilde{Z}$ on $^{k-1}N.$ Denote by $Z^{c^{k-1}}$
complete lift of a vector field $Z$ and by $\omega ^{v^{k-1}}$ vertical lift
of a 1-form $\omega $ defined on $N$ to $^{k-1}N.$ In (\ref{5.1}), if $%
\widetilde{Z}=Z^{c^{k-1}}$ and $\widetilde{\omega }=\omega ^{v^{k-1}},$ then
the \textit{vertical lift} of 1-form $\omega $ on $N$ to $^{k}N$ is the
1-form $\omega ^{v^{k}}$ on $^{k}N$ determined by
\begin{equation}
\omega ^{v^{k}}(Z^{c^{k}})=(\omega Z)^{v^{k}}.  \label{5.2}
\end{equation}%
\textbf{Proposition\thinspace 3.3.1} Let $\omega $ be a 1-form on $N$
locally expressed by (\ref{2.3}). Then $k$-th order vertical lift of $\omega
$ to $^{k}N$ is given by:
\begin{equation}
\omega ^{v^{k}}=dt+(\omega _{0i})^{v^{k}}dz^{0i}+(\overline{\omega }%
_{0i})^{v^{k}}d\overline{z}^{0i}.  \label{5.3}
\end{equation}%
\textbf{Proof: }Locally, we write $\omega ^{v^{k}}=dt+\omega _{ri}dz^{ri}+%
\overline{\omega }_{ri}d\overline{z}^{ri}$. Let $Z^{c^{k}}$ be complete lift
of a vector field $Z$ to extended manifold $^{k}N.$ Then from vertical lift
properties we get
\begin{eqnarray}
\omega ^{v^{k}}(Z^{c^{k}}) &=&(dt+\omega _{ri}dz^{ri}+\overline{\omega }%
_{ri}d\overline{z}^{ri})(Z^{c^{k}})  \label{5.3.1} \\
&=&1+\omega _{ri}\mathrm{C}_{r}^{k}(Z^{0i})^{v^{k-r}c^{r}}+\overline{\omega }%
_{ri}\mathrm{C}_{r}^{k}(Z^{0i})^{v^{k-r}c^{r}}  \notag
\end{eqnarray}%
and
\begin{eqnarray}
(\omega Z)^{v^{k}} &=&(dt+\omega _{0i}Z^{0i}+\overline{\omega }_{0i}%
\overline{Z}^{0i})^{v^{k}}  \label{5.3.2} \\
&=&1+(\omega _{0i})^{v^{k}}(Z^{0i})^{v^{k}}+(\overline{\omega }%
_{0i})^{v^{k}}(\overline{Z}^{0i})^{v^{k}}.  \notag
\end{eqnarray}%
By (\ref{5.2}), (\ref{5.3.1}), (\ref{5.3.2}) we have
\begin{equation*}
\omega _{ri}=0,\overline{\omega }_{ri}=0,\omega _{0i}=(\omega _{0i})^{v^{k}},%
\overline{\omega }_{0i}=(\overline{\omega }_{0i})^{v^{k}},1\leq r\leq k.
\end{equation*}%
Hence, the proof is complete.$\Box $

Let $^{k-1}N$ be $(k-1)$-th order extension of $N.$ Given by $\widetilde{%
\omega }$ a 1-form and by $\widetilde{Z}$ a vector field defined on $%
^{k-1}N. $ Then, the \textit{complete lift} of 1-form $\widetilde{\omega }$
on $^{k-1}N$ to $^{k}N$ is the 1-form $\widetilde{\omega }^{c}$ on $^{k}N$
defined by
\begin{equation}
\widetilde{\omega }^{c}(\widetilde{Z}^{c})=(\widetilde{\omega }\widetilde{Z}%
)^{c}.  \label{5.4}
\end{equation}%
By $Z^{c^{k-1}}$ and $\omega ^{c^{k-1}}$, let denote complete lifts of a
vector field $Z$ and a 1-form $\omega $ defined on $N$ to $^{k-1}N$ $.$ In (%
\ref{5.4}), if $\widetilde{Z}=Z^{c^{k-1}}$ and $\widetilde{\omega }=\omega
^{c^{k-1}},$ then the \textit{complete lift} of 1-form $\omega $ on $N$ to $%
^{k}N$ is the 1-form $\omega ^{c^{k}}$ on $^{k}N$ given by equality
\begin{equation}
\omega ^{c^{k}}(Z^{c^{k}})=(\omega Z)^{c^{k}}.  \label{5.5}
\end{equation}%
Similar to the proof of \textbf{Proposition\thinspace 3.3.1, }one can prove
the following:

\textbf{Proposition\thinspace 3.3.2} Let $\omega $ be a 1-form on $N$
locally expressed by (\ref{2.3})$.$ Then $k$-th order complete lift of $%
\omega $ to $^{k}N$ is given by:
\begin{equation}
\omega ^{c^{k}}=dt+(\omega _{0i})^{c^{k-r}v^{r}}dz^{ri}+(\overline{\omega }%
_{0i})^{c^{k-r}v^{r}}d\overline{z}^{ri}.  \label{5.6}
\end{equation}%
Let $Z$ be a vector field on manifold $N.$ Then the \textit{complete}-%
\textit{vertical lift of order }$(r,s)$ of $\omega \in \chi ^{\ast }(N)$ to $%
^{k}N$ is the 1-form $\omega ^{c^{r}v^{s}}\in \chi ^{\ast }(^{k}N)$ given by
equality
\begin{equation}
\omega ^{c^{r}v^{s}}(Z^{c^{k}})=(\omega Z)^{c^{r}v^{s}}.  \label{5.7}
\end{equation}%
\textbf{Proposition\thinspace 3.3.3} Let $N$ be a product manifold of
dimension $2m+1$ and $^{k}N$ its $k$-th order extension$.$ Suppose that the
1-form $\omega \in \chi ^{\ast }(N)$ is given by (\ref{2.3})$.$ Then
complete-vertical lift of order $(r,s)$ of $\omega $ to $^{k}N$ is
\begin{equation*}
\omega ^{c^{r}v^{s}}:\left( 1,\frac{\mathrm{C}_{l}^{r}}{\mathrm{C}_{l}^{k}}%
(\omega _{0i})^{v^{s+l}c^{r-l}},\frac{\mathrm{C}_{l}^{r}}{\mathrm{C}_{l}^{k}}%
(\overline{\omega }_{0i})^{v^{s+l}c^{r-l}}\right) ,0\leq l\leq k.
\end{equation*}%
\textbf{Proof: }Since $\omega ^{c^{r}v^{s}}$ is a 1-form on $^{k}N,$ with
respect to coordinate system $(t,z^{li},\,\overline{z}^{li})$ one writes $%
\omega ^{c^{r}v^{s}}=dt+\omega _{li}dz^{li}+\overline{\omega }_{li}d%
\overline{z}^{li}$. Let $Z^{c^{k}}$ be complete lift of order $k$ of vector
field $Z$ to $^{k}N.$ Then, from complete and vertical lift properties we
have
\begin{equation}
\omega ^{c^{r}v^{s}}(Z^{c^{k}})=\{1+\mathrm{C}_{l}^{r}\omega
_{li}(Z^{0i})^{v^{k-l}c^{l}}+\mathrm{C}_{l}^{r}\overline{\omega }_{li}(%
\overline{Z}^{0i})^{v^{k-l}c^{l}}\}  \label{5.8}
\end{equation}%
and
\begin{eqnarray}
(\omega Z)^{c^{r}v^{s}} &=&(1+\omega _{0i}Z^{0i}+\overline{\omega }_{0i}%
\overline{Z}^{0i})^{c^{r}v^{s}}  \label{5.9} \\
&=&1+\mathrm{C}_{h}^{r}(\omega _{0i})^{v^{s+h}c^{r-h}}(Z^{0i})^{v^{k-h}c^{h}}
\notag \\
&&+\mathrm{C}_{h}^{r}(\overline{\omega }_{0i})^{v^{s+h}c^{r-h}}(\overline{Z}%
^{0i})^{v^{k-h}c^{h}}.  \notag
\end{eqnarray}%
By (\ref{5.7}), (\ref{5.8}) and (\ref{5.9}), using $l=h$ from the following
equalities
\begin{equation*}
(Z^{0i})^{v^{k-l}c^{l}}=(Z^{0i})^{v^{k-h}c^{h}}\,{\ }\text{{and \ }}(%
\overline{Z}^{0i})^{v^{k-l}c^{l}}=(\overline{Z}^{0i})^{v^{k-h}c^{h}}
\end{equation*}%
we have
\begin{equation*}
\omega _{li}=(\frac{\mathrm{C}_{l}^{r}}{\mathrm{C}_{l}^{k}}(\omega
_{0i})^{v^{s+l}c^{r-l}},\overline{\omega }_{li}=(\frac{\mathrm{C}_{l}^{r}}{%
\mathrm{C}_{l}^{k}}(\overline{\omega }_{0i})^{v^{s+l}c^{r-l}},0\leq l\leq k
\end{equation*}%
Hence, the proof is complete.$\Box $

There exists the commutative property for complete-vertical lift of 1-forms.
Clearly, it means that complete-vertical lifts of order $(r,s)$ and
complete-vertical lifts of order $(s,r)$ of 1-forms on $N$ to its extended
manifold $^{k}N$ are equal. The property of complete-vertical lifts of order
$(r,s)$ of 1-forms on manifold $N$ is
\begin{equation*}
(f\omega )^{c^{r}v^{s}}=\mathrm{C}_{h}^{r}f^{v^{s+h}c^{r-h}}\omega
^{c^{h}v^{k-h}},\,0\leq r,s\leq k,\,(r+s=k).
\end{equation*}%
The \textit{horizontal lift} of a 1-form\textbf{\ }$\omega $ on $N$ to $%
^{k}N $ is the 1-form $\omega ^{H^{k}}$ on $^{k}N$ given by
\begin{equation*}
\omega ^{H^{k}}(Z^{H^{k}})=0,\,\omega ^{H^{k}}(Z^{v^{k}})=(\omega Z)^{v^{k}}.
\end{equation*}%
Considering $\omega =dt+\omega _{0i}dz^{0i}+\overline{\omega }_{0i}d%
\overline{z}^{0i},$ we obtain
\begin{equation*}
\omega ^{H^{k}}=dt+\omega _{ri}\eta ^{ri}+\overline{\omega }_{ri}\overline{%
\eta }^{ri},
\end{equation*}%
where $\eta ^{ri}=\overline{d}z^{r+1i}+\Gamma _{rj}^{ri}\overline{d}z^{ri},$
$\overline{\eta }^{ri}=\overline{d}\overline{z}^{r+1i}+\overline{\Gamma }%
_{rj}^{ri}\overline{d}\overline{z}^{ri},1\leq i,j\leq m.$ An \textit{%
extended coframe }adapted\textit{\ }to $\nabla $ on $^{k}N$ is the dual
coframe $\left\{ dt,\theta ^{ri}=dz^{ri},\overline{\theta }^{\alpha ri}=d%
\overline{z}^{ri},\eta ^{ri},\overline{\eta }^{ri}\right\} $.\thinspace The
properties of the higher order vertical, complete and horizontal lifts of
1-forms on $N$ are

$%
\begin{array}{ll}
i) &
\begin{array}{l}
(\omega +\lambda )^{v^{r}}=\omega ^{v^{r}}+\lambda ^{v^{r}},(\omega +\lambda
)^{c^{r}}=\omega ^{c^{r}}+\lambda ^{c^{r}}, \\
(\omega +\lambda )^{H^{k}}=\omega ^{H^{k}}+\lambda ^{H^{k}}%
\end{array}
\\
ii) &
\begin{array}{l}
(f\omega )^{v^{r}}=f^{v^{r}}\omega ^{v^{r}},(f\omega )^{c^{r}}=\sum_{j=0}^{r}%
\mathrm{C}_{j}^{r}f^{c^{r-j}v^{j}}\omega ^{c^{j}v^{r-j}}, \\
\omega ^{H^{k}}(Z^{H^{k}})=0,\,\omega ^{H^{k}}(Z^{v^{k}})=(\omega Z)^{v^{k}}%
\end{array}
\\
iii) & \chi ^{\ast }(V)=Sp\left\{ dt,dz^{0i},d\overline{z}^{0i}\right\} ,%
\text{ }\chi ^{\ast }(^{k}V)=Sp\left\{ dt,dz^{ri},d\overline{z}^{ri}\right\}
, \\
&
\begin{array}{l}
(dz^{0i})^{v^{r}}=dz^{0i},\,(d\overline{z}^{0i})^{v^{r}}=d\overline{z}%
^{0i},\,{\ }(dt)^{v^{r}}=dt, \\
(dz^{0i})^{c^{r}}=dz^{ri},\,(d\overline{z}^{0i})^{c^{r}}=d\overline{z}%
^{ri},\,(dt)^{c^{r}}=dt \\
(dt)^{H^{k}}=dt,(dz^{0i})^{H^{k}}=\eta ^{0i},\,(d\overline{z}^{0i})^{H^{k}}=%
\overline{\eta }^{0i},%
\end{array}%
\end{array}%
$

for all $\omega ,\lambda \in \chi ^{\ast }(N)$ and $f\in \mathcal{F}(N)$.

\subsection{\textbf{Higher Order Lifts of Tensor Fields of Type (1,1)}}

\textbf{\ }This subsection studies the extended definitions and properties
about higher order lifts of a tensor field of type (1,1) $\phi $ defined on $%
N$ to $^{k}N$.

Let $^{k-1}N$ be $(k-1)$-th order extension of $N$. Denote by $\widetilde{%
\phi }=\phi _{k-1}$ a tensor field of type (1,1), by $\widetilde{\xi }$ a
vector field and by $\widetilde{\eta }$ a 1-form defined on $^{k-1}N.$ Then
the \textit{vertical lift} of a tensor field of type (1,1) $\widetilde{\phi }
$ to $^{k}N$ is the tensor field $\widetilde{\phi }^{v}$ such that
\begin{equation}
\widetilde{\phi }^{v}(\widetilde{\xi }^{c})=(\widetilde{\phi }\widetilde{\xi
})^{v},\widetilde{\eta }^{v}(\widetilde{\phi }^{v})=(\widetilde{\eta }%
\widetilde{\phi })^{v}.  \label{6.1}
\end{equation}%
Now, let $Z^{c^{k-1}}$ and $\phi _{{}}^{v^{k-1}}$be respectively complete
lift of a vector field $Z\in \chi (N)$ and vertical lift of a tensor field
of type (1,1) $\phi \in \Im _{1}^{1}(N)$ to $^{k-1}N.$ In (\ref{6.1}), if $%
\widetilde{\xi }=\xi ^{c^{k-1}},\widetilde{\eta }=\eta ^{v^{k-1}}$ and $%
\widetilde{\phi }=\phi ^{v^{k-1}},$ then the \textit{vertical lift} of $\phi
$ to $^{k}N$ is the tensor field $\phi ^{v^{k}}$ on $^{k}N$ given by
\begin{equation}
\phi ^{v^{k}}(\xi ^{c^{k}})=(\phi \xi )^{v^{k}},\eta ^{v^{k}}(\phi
^{v^{k}})=(\eta \phi )^{v^{k}}.  \label{6.2}
\end{equation}%
Similarly, we define higher order complete lift of a tensor field of type
(1,1) $\phi $ on $N$. Denote by $\widetilde{\phi }$ a tensor field of type
(1,1), by $\widetilde{\xi }$ a vector field and by $\widetilde{\eta }$ a
1-form defined on $^{k-1}N.$ Then the \textit{complete lift} of $\widetilde{%
\phi }$ to $^{k}N$ is the structure $\widetilde{\phi }^{c}\in \Im
_{1}^{1}(^{k}N)$ determined by
\begin{equation}
\widetilde{\phi }^{c}(\widetilde{\xi }^{c})=(\widetilde{\phi }\widetilde{\xi
})^{c},\widetilde{\eta }^{c}(\widetilde{\phi }^{c})=(\widetilde{\eta }%
\widetilde{\phi })^{c}.  \label{6.3}
\end{equation}%
Presently, let $\xi ^{c^{k-1}},\eta ^{c^{k-1}}$ and $\phi ^{c^{k-1}}$be
respectively complete lifts of a vector field $\xi \in \chi (N),$a 1-form $%
\eta \in \chi (N)$ and a tensor field of type (1,1) $\phi \in \Im
_{1}^{1}(N) $ to $^{k-1}N.$ In (\ref{6.3}), if $\widetilde{\xi }=\xi
^{c^{k-1}},$ $\widetilde{\eta }=\eta ^{c^{k-1}}$ and $\widetilde{\phi }=\phi
^{c^{k-1}},$ then the \textit{complete lift }of $\phi \in \Im _{1}^{1}(N)$
to $^{k}N$ is the structure $\phi ^{c^{k}}$ on $^{k}N$ given by
\begin{equation}
\phi ^{c^{k}}(\xi ^{c^{k}})=(\phi \xi )^{v^{k}},\eta ^{c^{k}}(\phi
^{c^{k}})=(\eta \phi )^{c^{k}}.  \label{6.4}
\end{equation}

\end{document}